\documentclass[brochure,12pt]{bourbaki}

\date{Mars 2005}
\bbkannee{57\`eme ann\'ee, 2004-2005}
\bbknumero{947}
\title{Cat\'egories d\'eriv\'ees et g\'eom\'etrie birationnelle}
\subtitle{d'apr\`es Bondal, Orlov, Bridgeland, Kawamata...}
\author{Rapha\"el ROUQUIER}
\address{Institut de Math\'ematiques de Jussieu\\
UMR 7586 du CNRS\\
UFR de Math\'ematiques, Universit\'e Paris 7\\
Case 7012\\
2 place Jussieu\\
F--75251 PARIS C\'edex 05}
\email{rouquier@math.jussieu.fr}

\input xypic
\def\BA{{\mathbf{A}}}
\def\BC{{\mathbf{C}}}
\def\BG{{\mathbf{G}}}
\def\BL{{\mathbf{L}}}
\def\BP{{\mathbf{P}}}
\def\BQ{{\mathbf{Q}}}
\def\BR{{\mathbf{R}}}
\def\BZ{{\mathbf{Z}}}

\def\CA{{\mathcal{A}}}
\def\CB{{\mathcal{B}}}
\def\CD{{\mathcal{D}}}
\def\CE{{\mathcal{E}}}
\def\CF{{\mathcal{F}}}
\def\CG{{\mathcal{G}}}
\def\CH{{\mathcal{H}}}
\def\CI{{\mathcal{I}}}
\def\CJ{{\mathcal{J}}}
\def\CK{{\mathcal{K}}}
\def\CL{{\mathcal{L}}}
\def\CN{{\mathcal{N}}}
\def\CO{{\mathcal{O}}}
\def\CP{{\mathcal{P}}}
\def\CQ{{\mathcal{Q}}}
\def\CT{{\mathcal{T}}}
\def\CX{{\mathcal{X}}}

\def\Gm{{\mathfrak{m}}}

\def\tX{{\tilde{X}}}
\def\tY{{\tilde{Y}}}

\newcommand{\iso}{{\xrightarrow\sim}}

\def\ampl{\operatorname{ampl}\nolimits}

\def\Aut{\operatorname{Aut}\nolimits}

\def\ch{\operatorname{ch}\nolimits}
\def\codim{\operatorname{codim}\nolimits}
\def\coh{\operatorname{coh}\nolimits}
\def\mcoh{\operatorname{\!-coh}\nolimits}
\def\End{\operatorname{End}\nolimits}
\def\Ext{\operatorname{Ext}\nolimits}
\def\hd{\operatorname{hd}\nolimits}
\def\Hom{\operatorname{Hom}\nolimits}
\def\CHom{{{\mathcal H}om}}

\def\Per{\operatorname{Per}\nolimits}
\def\Pic{\operatorname{Pic}\nolimits}
\def\Proj{\operatorname{Proj}\nolimits}
\def\rang{\operatorname{rang}\nolimits}
\def\Spec{\operatorname{Spec}\nolimits}
\def\Supp{\operatorname{Supp}\nolimits}
\def\td{\operatorname{td}\nolimits}

\def\ie{{\em i.e.}}

\def\og{\guillemotleft}
\def\fg{\!\!\guillemotright}

\begin{document}
\maketitle
\section{Introduction}
L'objet principal de cet expos\'e est la cat\'egorie d\'eriv\'ee $D^b(X)$
des faisceaux coh\'erents sur une vari\'et\'e $X$.
La cat\'egorie d\'eriv\'ee organise l'information homologique
(groupes d'extensions entre faisceaux coh\'erents) et num\'erique
($K$-th\'eorie). Nous allons \'etudier son comportement au cours des 
op\'erations de \og chirurgie alg\'ebrique\fg\ (\og flips\fg\ et \og flops\fg).

La cat\'egorie d\'eriv\'ee d'un espace projectif se d\'ecrit \`a partir
d'une alg\`ebre de dimension finie (Beilinson, 1978) et ceci a plac\'e dans
un cadre appropri\'e les descriptions de fibr\'es vectoriels en terme d'alg\`ebre
lin\'eaire. \`A la suite de ce r\'esultat, des descriptions analogues 
(d\'ecomposition semi-orthogonale de la cat\'egorie d\'eriv\'ee)
ont \'et\'e recherch\'ees pour d'autres vari\'et\'es.
De telles d\'ecompositions devraient appara{\^\i}tre en pr\'esence d'un \og flip\fg,
\'etape cruciale du programme de Mori de mod\`eles minimaux (MMP)
pour la classification des vari\'et\'es projectives lisses,
et cela a amen\'e en particulier la question de l'invariance
de la cat\'egorie d\'eriv\'ee par \og flop\fg\ (Bondal-Orlov, 1995). 
D'un autre c\^ot\'e, la conjecture homologique de sym\'etrie miroir
(Kontsevich, 1994) a elle aussi pos\'e le probl\`eme de l'invariance
birationnelle de la cat\'egorie d\'eriv\'ee, pour des vari\'et\'es de Calabi-Yau.
Ind\'ependamment, la construction d'une \'equivalence d\'eriv\'ee entre une vari\'et\'e
ab\'elienne et sa duale (Mukai, 1981) a montr\'e la relation entre la r\'ealisation
d'une vari\'et\'e $X$ comme un espace de modules d'objets sur $Y$ et la
construction d'une \'equivalence (ou d'un foncteur pleinement fid\`ele)
entre la cat\'egorie d\'eriv\'ee de $X$ et celle de $Y$.

\medskip
Commen\c cons par poser des probl\`emes sur les cat\'egories
d\'eriv\'ees de vari\'et\'es, en suivant trois points importants du MMP.

On note $K_X$ le diviseur canonique d'une vari\'et\'e lisse $X$ et on
consid\`ere l'\'equivalence lin\'eaire entre diviseurs.
On se donne un diagramme o\`u $f$ et $g$ sont des morphismes
birationnels entre vari\'et\'es projectives lisses complexes
\begin{equation}
\label{Kequiv}
\xymatrix{
& Z\ar[dl]_f \ar[dr]^g \\
X && Y
}
\end{equation}

On a une premi\`ere conjecture sur les flops g\'en\'eralis\'es
(cf \cite[Conjecture 4.4]{BoOrICM} et \cite[Conjecture 1.2]{KaDK}):
\begin{conj}[Bondal-Orlov]
\label{conjflop}
Si $f^*K_X\sim g^*K_Y$,
alors $D^b(X)\simeq D^b(Y)$.
\end{conj}

On sait que dans cette situation les nombres de Hodge
co\"{\i}ncident (cf remarque \ref{Hodge}).
La conjecture a une r\'eponse positive
en dimension $3$ (corollaire \ref{CY3}), pour
des vari\'et\'es symplectiques de dimension $4$ (corollaire \ref{sym4}) et
dans le cadre torique (th\'eor\`eme \ref{torique}).
Pour des vari\'et\'es de Calabi-Yau ($\omega$ trivial), on s'attend donc
\`a ce que birationalit\'e et \'equivalence d\'eriv\'ee co{\"\i}ncident,
comme le pr\'edit la conjecture de Kontsevich de sym\'etrie miroir \cite{Ko}.
La \og r\'eciproque\fg\ de la conjecture \ref{conjflop} n'est pas vraie
(remarque \ref{Uehara}).

\smallskip
On a une conjecture sur les flips g\'en\'eralis\'es \cite[Conjecture 4.4]{BoOrICM}:
\begin{conj}[Bondal-Orlov]
\label{conjflip}
Si $f^*K_X-g^*K_Y$ est \'equivalent \`a un diviseur effectif, alors il existe un
foncteur pleinement fid\`ele $D^b(Y)\to D^b(X)$.
\end{conj}

La minimisation d'une vari\'et\'e dans le MMP devrait alors s'interpr\'eter comme une
minimisation de la cat\'egorie d\'eriv\'ee, un mod\`ele minimal pour une vari\'et\'e $X$
devant \^etre construit comme
un espace de module d'objets de la cat\'egorie d\'eriv\'ee de $X$.
Dans la conjecture \ref{conjflip}, il
serait aussi souhaitable de savoir d\'ecrire l'orthogonal de l'image 
de $D^b(Y)$ dans $D^b(X)$.

\smallskip
On a enfin une conjecture de finitude \cite[Conjecture 1.5]{KaDK}:
\begin{conj}[Kawamata]
Soit $X$ une vari\'et\'e projective lisse. Alors, il n'existe qu'un nombre
fini de classes d'isomorphisme de vari\'et\'es projectives lisses
$Y$ telles que $D^b(X)\simeq D^b(Y)$.
\end{conj}

La r\'eponse est positive en dimension $\le 2$ (cf \cite[Corollary 1.2]{BrMa}
et \cite[Theorem 1.6]{KaDK}) et pour $X,Y$ des vari\'et\'es ab\'eliennes
\cite[Corollary 2.8]{Orab}.

\smallskip
Dans la premi\`ere partie de cet expos\'e, nous montrons,
dans la situation extr\^eme o\`u le fibr\'e canonique est ample ou anti-ample,
comment reconstruire la vari\'et\'e \`a partir de la cat\'egorie d\'eriv\'ee
(th\'eor\`eme \ref{BO}).
Nous expliquons ensuite le m\'ecanisme de d\'evissage de cat\'egories d\'eriv\'ees
(\S \ref{secsemi}). Dans le
\S \ref{secFM}, nous pr\'esentons la
construction de transformations \`a noyau et les invariants transport\'es
par les \'equivalences, puis nous exposons le cas des vari\'et\'es ab\'eliennes
et des surfaces, o\`u la th\'eorie est presque compl\`ete. La partie
\S \ref{flipflop} pr\'esente plusieurs cas de r\'eponse positive aux
conjectures \ref{conjflop} et \ref{conjflip}.

\smallskip
Outre le livre en pr\'eparation \cite{Huy}, le lecteur pourra consulter
\cite{BoOrICM,Ca3,HiVdB,Or2} pour des expos\'es g\'en\'eraux.

\smallskip
{\em Je remercie Arnaud Beauville, Tom Bridgeland, Olivier Debarre et
Alastair King pour leurs remarques sur une version pr\'eliminaire de ce texte,
et Daniel Huybrechts, pour de très nombreuses discussions.}

\section{Propri\'et\'es internes}

\subsection{Reconstruction}
\subsubsection{Terminologie}
Une {\em vari\'et\'e} sera pour nous un sch\'ema quasi-projectif $X$
sur $\BC$. La plupart du temps, il s'agira de vari\'et\'es projectives lisses.
Les points consid\'er\'es seront toujours des points ferm\'es.

La cat\'egorie d\'eriv\'ee $D^b(X)$ est d\'efinie comme la localisation de la cat\'egorie
des complexes born\'es de faisceaux coh\'erents en la classe des
quasi-isomorphismes (=morphismes de complexes
qui induisent un isomorphisme entre faisceaux de cohomologie).
Ses objets sont donc les complexes born\'es de faisceaux
coh\'erents sur $X$. Les fl\`eches sont obtenues \`a partir de morphismes
de complexes auxquels les inverses des quasi-isomorphismes ont
\'et\'e ajout\'es.
La cat\'egorie $D^b(X)$ n'est pas ab\'elienne en g\'en\'eral, mais elle
poss\`ede la structure de cat\'egorie triangul\'ee~: le r\^ole des
suites exactes courtes est jou\'e par les triangles distingu\'es.

Tous les foncteurs consid\'er\'es entre cat\'egories triangul\'ees seront
triangul\'es. Une sous-cat\'egorie {\em \'epaisse} d'une cat\'egorie
triangul\'ee est une sous-cat\'egorie triangul\'ee pleine close par facteurs
directs. La sous-cat\'egorie d'une cat\'egorie triangul\'ee {\em engendr\'ee}
(resp. {\em faiblement engendr\'ee})
par une famille d'objets est la plus petite sous-cat\'egorie pleine
triangul\'ee (resp. la plus petite sous-cat\'egorie \'epaisse)
contenant cette famille.

\subsubsection{Cat\'egories ab\'eliennes}
Soit $X$ une vari\'et\'e lisse.
On sait, depuis Gabriel \cite{Ga}, que la cat\'egorie des faisceaux coh\'erents
$X\mcoh$ sur $X$ d\'etermine $X$:
\begin{itemize}
\item l'application qui a un ferm\'e associe les faisceaux support\'es par ce
ferm\'e induit une bijection $Z\mapsto \CI_Z$ de l'ensemble des ferm\'es de $X$
vers l'ensemble des sous-cat\'egories de Serre de $X\mcoh$ engendr\'ees
par un objet
\item
la cat\'egorie quotient $X\mcoh/\CI_Z$ est \'equivalente \`a $(X-Z)\mcoh$
et son centre s'identifie \`a $\CO_X(X-Z)$.
\end{itemize}

Suivant Thomason et Balmer \cite{Ba} (voir aussi \cite[Theorem 3.11]{Rou}), ce 
principe de reconstruction s'\'etend \`a la cat\'egorie $D^b(X)$, si
on munit celle-ci de sa structure tensorielle, en plus de sa structure
triangul\'ee:
\begin{itemize}
\item
l'application qui \`a un ferm\'e $Z$ de $X$ associe la sous-cat\'egorie pleine
$D_Z^b(X)$ de $D^b(X)$ des complexes dont les faisceaux de cohomologie
sont support\'es par $Z$ est injective, d'image l'ensemble
des sous-cat\'egories \'epaisses faiblement
engendr\'ees par un \'el\'ement et $\otimes$-id\'eales (\ie,
stables par $-\otimes L$ pour tout $L\in D^b(X)$)
 \cite[Theorem 3.15]{Thclassif}
\item La cat\'egorie
$D^b(X-Z)$ s'identifie \`a $D^b(X)/D^b_Z(X)$ et, si $X-Z$ est un
ouvert affine, le centre de $D^b(X-Z)$ s'identifie \`a $\CO_X(X-Z)$.
\end{itemize}

Par cons\'equent, $X\mcoh$, vue comme cat\'egorie ab\'elienne, et $D^b(X)$,
vue comme cat\'egorie triangul\'ee tensorielle, ne sont pas des invariants
int\'eressants de $X$ !

\subsubsection{Fibr\'e canonique (anti-)ample}
La cat\'egorie $D^b(X)$, munie de sa seule structure de cat\'egorie
triangul\'ee, ne d\'etermine pas la vari\'et\'e $X$, mais elle appara{\^\i}t comme
un invariant int\'eressant. Dans la suite, les cat\'egories d\'eriv\'ees seront
consid\'er\'ees avec leur seule structure triangul\'ee (voir \`a ce sujet
\S \ref{axiomes}).
Le premier exemple (Mukai) est celui de
l'\'equivalence d\'eriv\'ee entre une vari\'et\'e ab\'elienne et sa duale (th\'eor\`eme
\ref{equivabelienne}).

Bondal et Orlov d\'emontrent que lorsque le fibr\'e canonique
$\omega_X$ est ample ou anti-ample, alors la vari\'et\'e est
d\'etermin\'ee par sa cat\'egorie d\'eriv\'ee \cite[Theorem 2.5]{BoOr2}.

\begin{theo}[Bondal-Orlov]
\label{BO}
Soit $X$ une vari\'et\'e projective lisse telle que
$\omega_X$ ou $\omega_X^{-1}$ est ample.
Si $Y$ est une vari\'et\'e projective
lisse et si on a une \'equivalence de cat\'egories triangul\'ees
$D^b(X\mcoh)\simeq D^b(Y\mcoh)$, alors $X\simeq Y$.
\end{theo}

Un point crucial est jou\'e dans la preuve par la notion de foncteur de
Serre \cite[\S 3]{BoKa}.
Soit $\CT$ une cat\'egorie $\BC$-lin\'eaire telle que
$\dim\Hom(M,N)<\infty$ pour tous $M,N\in\CT$. Un foncteur de Serre
est la donn\'ee d'une auto-\'equivalence $S:\CT\iso\CT$ et d'isomorphismes
bifonctoriels pour tous $M,N\in\CT$:
$$\Hom(M,N)\iso\Hom(N,SM)^*.$$
Si on voit $\CT$ comme une \og alg\`ebre avec plusieurs objets\fg, alors
ceci correspond \`a la notion d'alg\`ebre de Frobenius.

Un foncteur de Serre, s'il existe, est unique \`a isomorphisme unique pr\`es.
Par cons\'equent, une \'equivalence de cat\'egories commute avec les
foncteurs de Serre. En outre, si
$\CT$ est une cat\'egorie triangul\'ee, alors un foncteur de Serre est
automatiquement triangul\'e.

\smallskip
La d\'efinition est motiv\'ee par la dualit\'e de Serre:
\begin{theo}
Si $X$ est une vari\'et\'e projective lisse purement de dimension $n$, alors
$S=\omega_X[n]\otimes -$ est un foncteur de Serre pour $D^b(X)$.
\end{theo}

\smallskip
\noindent{\sc Preuves du th\'eor\`eme \ref{BO}} (esquisses) ---

On se ram\`ene facilement au cas o\`u $X$ et $Y$ sont connexes
(cf proposition \ref{indec})
et on fixe une \'equivalence $F:D^b(X)\iso D^b(Y)$. On suppose
$\omega_X$ ample (preuve identique dans l'autre cas).
Le th\'eor\`eme r\'esulterait imm\'ediatement de l'invariance des alg\`ebres
canoniques, si on savait que $\omega_Y$ \'etait ample
(cf \S \ref{invariants}).

$\bullet\ $1\`ere approche \cite{BoOr2}.
Sur une vari\'et\'e projective lisse connexe
$Z$, pour tout point $z$ et tout $i\in\BZ$,
les $C=\CO_z[i]\in D^b(Z)$\footnote{$\CO_z$ est le faisceau gratte-ciel
en $z$} v\'erifient
\begin{equation}
\label{objetpoint}
S(C)\simeq C[\dim Z],\ \End_{D^b(Z)}(C)=\BC \textrm{ et }
\Hom(C,C[i])=0\textrm{ pour }i<0.
\end{equation}
Sur la vari\'et\'e $X$ o\`u $\omega_X$ est ample, les conditions
(\ref{objetpoint}) caract\'erisent les objets $\CO_x[i]$ dans $D^b(X)$.
On en d\'eduit que 
l'ensemble $\{F(\CO_x)[i]\}_{i,x}$ contient les $\CO_y[j]$ pour $y\in Y$ et
$j\in\BZ$. Si $F(\CO_x)$ n'est pas de cette forme, il est
orthogonal aux $\CO_y[j]$, donc il est nul. On d\'eduit alors que
$F$ envoie tout $\CO_x[i]$ sur un $\CO_y[j]$ et ceci
induit une bijection entre points de $X$ et de $Y$.
On caract\'erise ensuite les faisceaux inversibles d\'ecal\'es sur une vari\'et\'e
lisse $Z$ comme les $C\in D^b(Z)$ tels que pour tout $z\in Z$, il existe
$n\in\BZ$ tel que
$$\Hom(L,\CO_z[n])\simeq\BC\textrm{ et }
\Hom(L,\CO_z[i])=0 \textrm{ pour } i\not=n.$$
On en d\'eduit que $F$ envoie un faisceau inversible sur un faisceau inversible
d\'ecal\'e. Soit $L\in\Pic(X)$. Quitte \`a d\'ecaler $F$, alors on peut
supposer $F(L)\in\Pic(Y)$.
L'alg\`ebre $\bigoplus_{i\ge 0}\Hom(L,S^i(L)[-i\dim X])$ est isomorphe
\`a l'alg\`ebre canonique de $X$ et les ouverts d\'efinis par ses \'el\'ements
forment
une base de la topologie de $X$. Cette alg\`ebre est isomorphe \`a l'alg\`ebre
d\'efinie de la m\^eme fa\c con pour $Y$ et elle donne donc une base de la
topologie de $Y$. Ceci montre que $\omega_Y$ est ample et que
les alg\`ebres canoniques de $X$ et $Y$ sont isomorphes.

\smallskip
$\bullet\ $2\`eme approche \cite[\S 4]{HiVdB}.
On commence comme ci-dessus par v\'erifier que les $\CO_x[i]$
s'envoient sur des $\CO_y[j]$.
La suite de la preuve n'utilise plus que $\omega_X$ est ample.
On utilise le th\'eor\`eme \ref{repOrlov} plus bas qui affirme
qu'il existe $K\in D^b(Y\times X)$ tel que $F=\Phi_K$.
Alors, le lemme \ref{FMtrivial} plus bas montre que $Y\simeq X$.

\smallskip
$\bullet\ $3\`eme approche \cite[\S 3.2.4]{Rou}.
Soit $\CI$ une sous-cat\'egorie \'epaisse de $D^b(X)$. Si
$\CI$ est stable par $L^{-1}\otimes -$ pour un faisceau ample $L$, alors elle
est $\otimes$-id\'eale. Cette propri\'et\'e
est donc \'equivalente \`a la stabilit\'e sous $S^{-1}$. Par
cons\'equent, l'ensemble des ferm\'es de $X$ se retrouve \`a partir de $D^b(X)$
(\`a partir de sa seule structure triangul\'ee). 
Pour tout ferm\'e $Z$ de $Y$, il existe donc un ferm\'e $Z'$ de $X$ tel que
$F(D^b_{Z'}(X))=D_Z^b(Y)$. On montre que cette injection de l'ensemble
des ferm\'es de $Y$ vers ceux de $X$ se restreint en une bijection
$Y\to X$ d'inverse continu.
On identifie enfin les faisceaux d'anneaux.$\Box$

\smallskip
Bondal et Orlov \cite[Theorem 3.1]{BoOr2} d\'eterminent le 
groupe $\Aut(D^b(X))$ des classes d'isomorphisme
d'auto-\'equivalences de $D^b(X)$ lorsque $\omega_X^\pm$ est ample
(ceci est par exemple fourni par la deuxi\`eme preuve du th\'eor\`eme \ref{BO})~:
$$\Aut(D^b(X))=\Pic(X)\rtimes\Aut(X)\times\BZ.$$

\subsection{D\'ecompositions semi-orthogonales}
\label{secsemi}
\subsubsection{D\'ecompositions partielles}
Consid\'erons la forme d'Euler sur la $K$-th\'eorie $K_0(X)$
$$\langle [\CF],[\CG]\rangle=\sum_{i\ge 0}(-1)^i\dim\Ext^i(\CF,\CG)$$
Nous allons d\'ecrire l'analogue, pour la cat\'egorie d\'eriv\'ee, d'une base
triangulaire pour cette forme, o\`u plus g\'en\'eralement d'une d\'ecomposition
semi-orthogonale de $K_0(X)$.

\smallskip
Soit $\CI$ une sous-cat\'egorie \'epaisse
d'une cat\'egorie triangul\'ee $\CT$. On pose
${^\perp\CI}=\{C\in \CT| \Hom(C,I)=0\textrm{ pour tout }I\in\CI\}$
et
$\CI^\perp=\{C\in \CT| \Hom(I,C)=0\textrm{ pour tout }I\in\CI\}$.
On dit
que $\langle \CI^\perp,\CI\rangle$ est une {\em d\'ecomposition
semi-orthogonale} de $\CT$ lorsque pour tout objet $C$ de $\CT$, il existe
un triangle distingu\'e $C_1\to C\to C_2\rightsquigarrow$ avec
$C_1\in \CI$ et $C_2\in\CI^\perp$.
Ceci revient \`a demander que le foncteur canonique
$\CI^\perp\to \CT/\CI$ soit une \'equivalence ou \`a demander que
le foncteur d'inclusion $\CI\to\CT$ ait un adjoint \`a droite.

Lorsque $\CI^\perp=\langle \CK,\CJ\rangle$, on \'ecrit
$\CT=\langle \CK,\CJ,\CI\rangle$ et on g\'en\'eralise aux d\'ecompositions
$\CT=\langle\CI_1,\ldots,\CI_m\rangle$.

L'existence d'un \og g\'en\'erateur fort\fg\ pour $D^b(X)$ fournit un
th\'eor\`eme de repr\'esentabilit\'e \`a la Brown pour les foncteurs cohomologiques sur
$D^b(X)$ (cf \cite{BoKa} et \cite{BoVdB}) et on obtient un th\'eor\`eme
g\'en\'eral d'existence de d\'ecompositions:
\begin{theo}[Bondal, Kapranov, Van den Bergh]
\label{admissible}
Soient $X$ une vari\'et\'e projective lisse et
$\CI=D^b(X)$ une sous-cat\'egorie triangul\'ee pleine d'une cat\'egorie
triangul\'ee $\CT$. Alors,
le foncteur d'inclusion $\CI\to\CT$ a des adjoints \`a gauche et \`a droite,
\ie,
on a des d\'ecompositions semi-orthogonales
$\CT=\langle \CI^\perp,\CI\rangle$ et
$\CT=\langle \CI,{^\perp\CI}\rangle$.
\end{theo}

Une d\'ecomposition orthogonale de $D^b(X)$ correspond \`a une d\'ecomposition
de $X$ en union de composantes connexes:
\begin{prop}
\label{indec}
Soit $X$ une vari\'et\'e connexe. Soient $\CI_1$ et $\CI_2$ deux sous-cat\'egories
\'epaisses de $D^b(X)$ telles que $D^b(X)=\CI_1\oplus\CI_2$ (\ie,
$D^b(X)=\langle\CI_1,\CI_2\rangle=\langle\CI_2,\CI_1\rangle$).
Alors, $\CI_1=0$ ou $\CI_2=0$.
\end{prop}

\noindent{\sc Preuve} ---
Un objet ind\'ecomposable de $D^b(X)$ est dans $\CI_1$ ou dans $\CI_2$. Soient
$r,s$ tels que $\CO_X\in\CI_r$ et $\{r,s\}=\{1,2\}$.
Soit $X_i=\{x\in X\ |\ \CO_x\in\CI_i\}$.
Si $x\in X_s$, alors $\Hom(\CO_X,\CO_x)=0$, ce qui est impossible,
donc $X_s=\emptyset$. Si $C\in \CI_s$, alors
$\Hom(C,\CO_x[i])=0$ pour tout $x\in X$ et tout $i\in\BZ$, donc
$C=0$.$\Box$

\begin{rema}
Soit $\CT$ une cat\'egorie triangul\'ee avec un foncteur de Serre $S$ et
$\CT=\langle\CI,\CI^\perp\rangle$ une d\'ecomposition semi-orthogonale.
Alors, $\CT=\langle {^\perp\CI},\CI\rangle$ et
${^\perp\CI}=S^{-1}(\CI^\perp)$.

Consid\'erons en particulier $\CT=D^b(X)$ o\`u $X$ une vari\'et\'e projective lisse
connexe de Calabi-Yau. Alors,
il n'y a pas de d\'ecomposition semi-orthogonale non triviale de $D^b(X)$,
car une telle d\'ecomposition serait une d\'ecomposition orthogonale.
\end{rema}

\subsubsection{D\'ecompositions compl\`etes}
Voyons le cas particulier d'une {\em suite exceptionnelle} d'objets.
C'est une suite $(C_1,\ldots,C_m)$ d'objets de $\CT$
telle que
\begin{itemize}
\item
$\Hom(C_i,C_j[r])=0$ pour $r\in\BZ$ et $i>j$
\item
$\Hom(C_i,C_i[r])=0$ pour $r\not=0$
\item
$\End(C_i)=\BC$.
\end{itemize}

On dit que la suite est {\em compl\`ete} si elle engendre $\CT$.

\smallskip
Soit $(C_1,\ldots,C_m)$ une suite exceptionnelle compl\`ete.
Notons $\CI_i$ la sous-cat\'egorie triangul\'ee de $\CT$ engendr\'ee par
$C_i$. Nous noterons $D^b(A)$ la cat\'egorie d\'eriv\'ee born\'ee des
modules de type fini sur une alg\`ebre $A$.

On a une \'equivalence $C_i\otimes_{\BC} -:
D^b(\BC)\iso \CI_i$ et une d\'ecomposition semi-orthogonale
$\CT=\langle \CI_1,\ldots,\CI_m\rangle$.
R\'eciproquement, toute d\'ecomposition semi-orthogonale en des cat\'egories
\'equivalentes \`a $D^b(\BC)$ provient d'une suite exceptionnelle
d'objets.
L'ensemble $\{[C_i]\}$ forme une base de $K_0(\CT)$.
Si $\CT$ est localement de type fini (\ie, si
$\dim\bigoplus_i \Hom(M,N[i])<\infty$ pour tous $M,N\in\CT$), alors
la matrice de la forme d'Euler dans la base $\{[C_i]\}$ est triangulaire.

\begin{exem}
\label{Beilinson}
Soit $X=\BP^n$. Beilinson \cite{Be} montre que
$(\CO(-n),\CO(-n+1),\ldots,\CO(0))$ est une suite exceptionnelle compl\`ete.
L'orthogonalit\'e est claire. 
La r\'esolution de la diagonale $\Delta\subset \BP^n\times\BP^n$:
$$0\to \CO(-n)\boxtimes \Omega^n(n)\to\cdots\to\CO(-1)\boxtimes\Omega^1(1)
\to\CO\boxtimes\CO\to \CO_{\Delta}\to 0$$
d\'ecrit le foncteur identit\'e de $D^b(\BP^n)$ comme extension de foncteurs
$\CO(-i)\otimes H^*(\Omega^i(i)\otimes -)$ et ceci d\'emontre l'engendrement.
Soient $\CF=\bigoplus_{i=0}^n \CO(-i)$ et $A=\End(\CF)$, une alg\`ebre de
dimension finie. Alors, on a en plus ici
$\Ext^{>0}(\CF,\CF)=0$ et on d\'eduit que le foncteur
$R\Hom(\CF,-):D^b(\BP^n)\to D^b(A)$ est une \'equivalence.
\end{exem}

\subsubsection{Minimisation}
Kapranov a construit des suites exceptionnelles compl\`etes pour les
quadriques projectives lisses et les vari\'et\'es de drapeaux de type $A$
\cite{Kap} (cf \cite{Ru} pour un survol des constructions pour les
vari\'et\'es de Fano).

King \cite[Conjecture 9.3]{Ki}
conjecture que toute vari\'et\'e torique compl\`ete lisse $X$
admet une suite $(L_1,\ldots,L_n)$ de fibr\'es en droites telle que
$\Ext^{>0}(L_i,L_j)=0$ pour tous $i,j$ et les $L_i$ engendrent $D^b(X)$
(alors, $D^b(X)$ est \'equivalente \`a $D^b(A)$, o\`u
$A=\End(\bigoplus_i L_i)$). Kawamata \cite{Katoric} d\'emontre l'existence
d'une suite exceptionnelle compl\`ete de faisceaux pour toute vari\'et\'e
torique projective lisse.

D\`es que $K_0(X)$ n'est pas de type fini, $D^b(X)$ ne peut \^etre
\'equivalente \`a la cat\'egorie d\'eriv\'ee d'une alg\`ebre de dimension finie.
On montre par contre que pour toute vari\'et\'e $X$, il existe une dg-alg\`ebre
(=alg\`ebre diff\'erentielle gradu\'ee)
$A$ dont la cat\'egorie des complexes parfaits est \'equivalente \`a $D^b(X)$
(cf Keller, Thomason, Neeman, Kontsevich, Bondal-Van den Bergh,
\cite[Proposition 3.14 et Theorem 7.39]{Rou2}).

La recherche de mod\`eles minimaux est reli\'ee \`a la minimisation de la cat\'egorie
d\'eriv\'ee. On cherche une suite exceptionnelle
$L_1,\ldots,L_n$ avec $n$ maximal. Soit $\CI$ la sous-cat\'egorie
triangul\'ee de $D^b(X)$ engendr\'ee par les $L_i$. Alors, 
$D^b(X)=\langle \CI,\CI^\perp\rangle$ et la g\'eom\'etrie de $X$ devrait \^etre
en partie contr\^ol\'ee par la cat\'egorie triangul\'ee $\CT=\CI^\perp$. 
Il serait int\'eressant d'\'etudier l'ind\'ependance de $\CT$ du
choix de $L_1,\ldots,L_n$ et m\^eme son ind\'ependance birationnelle
(cf \cite{KnV14}).
La cat\'egorie $\CT$ appara{\^\i}t parfois comme la cat\'egorie d\'eriv\'ee d'une vari\'et\'e
$X'$, de dimension inf\'erieure ou \'egale
(cf \cite{KnV12} pour un exemple de vari\'et\'e de Fano $X$ de
dimension $3$ o\`u $X'$ est une courbe de genre $7$).
L'exemple le plus simple
est celui d'un \'eclatement de centre un espace projectif
(cf th\'eor\`eme \ref{eclatement}).

\section{Comparaisons}
\label{secFM}
\subsection{Transformations \`a noyau}
\subsubsection{D\'efinition}
L'id\'ee des transformations \`a noyau est la suivante: on se donne
une fonction $\phi:X\times Y\to \BC$. On a alors une application
des fonctions sur $Y$ vers les fonctions sur $X$ donn\'ee
par $f\mapsto (x\mapsto \int_Y f(y)\phi(x,y)dy)$.

Cette construction a un analogue pour les faisceaux coh\'erents
(les m\^emes constructions pour
les faisceaux constructibles ou les $\CD$-modules sont classiques).
Soient $X$ et $Y$ deux vari\'et\'es projectives lisses et 
$p:X\times Y\to X$, $q:X\times Y\to Y$ les deux projections.
$$
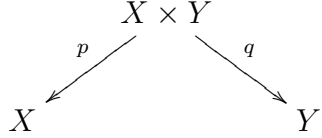
\xymatrix{
& X\times Y \ar[dl]_p \ar[dr]^q \\
X && Y
}$$

Soit $K\in D^b(X\times Y)$.
On d\'efinit alors le foncteur (dit de \og Fourier-Mukai\fg)
$\Phi_K:D^b(Y)\to D^b(X)$ par
$$\Phi_K(C)=\BR p_* (K\otimes^\BL q^*C).$$
Soit $K^\vee=\BR \CHom(K,\CO_{X\times Y})\in D^b(Y\times X)$.
Si $X$ et $Y$ sont de dimension pure, alors
les foncteurs $\Phi_{K^\vee\otimes p^*\omega_X[\dim X]}$
et $\Phi_{K^\vee\otimes q^*\omega_Y[\dim Y]}$ sont respectivement adjoints \`a
gauche et \`a droite de $\Phi_K$.

Soient $Z$ une autre vari\'et\'e projective lisse et $L\in D^b(Y\times Z)$.
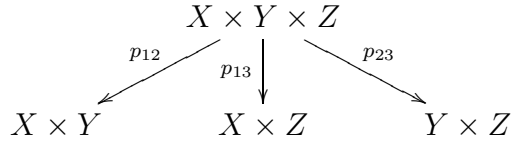
$$\xymatrix{
& X\times Y\times Z\ar[dl]_{p_{12}} \ar[dr]^{p_{23}}\ar[d]_{p_{13}}\\
X\times Y & X\times Z& Y\times Z
}$$
On pose $K\circ L=\BR p_{13*}(p_{12}^*K\otimes^\BL p_{23}^*L)$.
On a alors un isomorphisme canonique
$\Phi_K\circ\Phi_L\iso \Phi_{K\circ L}$.

\smallskip
Le lemme classique suivant permet de reconna{\^\i}tre quand $K$ provient d'un
isomorphisme de vari\'et\'es (cf \cite[Corollary 4.3]{HiVdB}).

\begin{lemm}
\label{FMtrivial}
On suppose que $Y$ est connexe et que pour tout $y\in Y$, il existe
$x\in X$ et $n\in\BZ$ tels que $\Phi_K(\CO_y)\simeq \CO_x[n]$.
Alors, il existe un morphisme $\sigma:Y\to X$
de graphe $\Gamma_\sigma$ et il existe $L\in\Pic(Y)$ et $m\in\BZ$ tels
que $K\simeq \CO_{\Gamma_\sigma}\otimes q^*L[m]$.

Si $\Phi_K$ est une \'equivalence, alors $\sigma$ est un isomorphisme.
\end{lemm}

\noindent{\sc Preuve} (esquisse) ---
Soit $y\in Y$ d'anneau local $\CO_{\Gm_y}$ et soit $n\in\BZ$
tel que $K\otimes_{\CO_Y}^\BL\CO_y$ est concentr\'e en degr\'e $-n$.
Le lemme de Nakayama
montre que $q_*(K\otimes_{\CO_Y}\CO_{\Gm_y})\simeq\CO_{\Gm_y}[n]$.
Il existe alors un voisinage ouvert $U$ de $y$ tel que
$q_*(K\otimes_{\CO_Y}\CO_U)\simeq \CO_U[n]$ et on obtient
un morphisme $U\to X$. Ceux-ci se recollent en $\sigma:Y\to X$ avec
les propri\'et\'es voulues.
Si $\Phi_K$ est une \'equivalence, alors 
$K^\vee\otimes p^*\omega_X[\dim X]$ d\'efinit un morphisme $X\to Y$
inverse de $\sigma$.$\Box$

\subsubsection{Pleine fid\'elit\'e}
\label{secplf}
La pleine fid\'elit\'e du foncteur $\Phi_K$ peut se tester sur
une famille d'objets appropri\'ee (cf \cite[Theorems 5.1 and 5.4]{Br1}).
Nous utilisons ici les faisceaux gratte-ciel \cite[Theorem 1.1]{BoOr1},
le point-clef \'etant qu'un objet orthogonal (\`a gauche ou
\`a droite) aux faisceaux gratte-ciel et \`a leurs d\'ecal\'es est nul.

\begin{prop}[Bondal-Orlov]
\label{critere}
Soit $K\in D^b(X\times Y)$.
Le foncteur $\Phi_K$ est pleinement fid\`ele si et seulement si
pour tous $y,y'\in Y$, on a 
$$\Hom(\Phi_K(\CO_y),\Phi_K(\CO_{y'})[i])=
\begin{cases}
0 & \text{ sauf si } y=y' \text{ et } 0\le i\le \dim Y \\
k & \text{ si } y=y' \text{ et }i=0.
\end{cases}$$
C'est une \'equivalence si en plus
$\Phi_K(\CO_y)\otimes\omega_X\simeq \Phi_K(\CO_y)$
pour tout $y\in Y$.
\end{prop}

Cette proposition montre que lorsque $X$ et $Y$ ont des ferm\'es stricts $Z$
et $Z'$ en dehors desquels $K$ est le faisceau de structure du
graphe d'un isomorphisme
$Y-Z'\iso X-Z$, alors le crit\`ere pr\'ec\'edent peut se v\'erifier en
rempla\c cant $X$ et $Y$ par leurs compl\'et\'es formels le long de $Z$ et $Z'$.
Ceci permet de substituer \`a $X$ et $Y$ des mod\`eles pr\'ef\'er\'es, \`a condition
de garder les m\^emes compl\'et\'es formels (cf la preuve du th\'eor\`eme
\ref{flopMukai}).

\smallskip
Un cas particulier de pleine fid\'elit\'e est fourni par le
r\'esultat suivant, qui se d\'eduit imm\'ediatement de la formule de projection.

\begin{prop}
\label{pleinefideliteprojection}
Soit $f:V\to W$ un morphisme entre vari\'et\'es projectives lisses.
Si le morphisme canonique $\CO_W\to\BR f_*\CO_V$ est un isomorphisme, alors
$\BL f^*:D^b(W)\to D^b(V)$ est pleinement fid\`ele.
\end{prop}

\subsubsection{Familles}

Soient $X',Y'$ deux vari\'et\'es projectives lisses et $K'\in D^b(X'\times Y')$.
Alors, on a le r\'esultat classique (cf \cite[Proposition 2.1.7]{Or2}):
\begin{prop}
Si $\Phi_K:D^b(Y)\to D^b(X)$ et 
$\Phi_{K'}:D^b(Y')\to D^b(X')$ sont pleinement fid\`eles (resp. sont des
\'equivalences), alors $\Phi_{K\boxtimes K'}:D^b(Y\times Y')\to D^b(X\times X')$
est pleinement fid\`ele (resp. est une \'equivalence).
\end{prop}

Soient $p:X\to S$ et $q:Y\to S$ des morphismes projectifs lisses entre
vari\'et\'es projectives lisses. Soient $s_0\in S$, $X_0=p^{-1}(s_0)$ et
$Y_0=q^{-1}(s_0)$. Soient $i:X_0\to X$, $j:Y_0\to Y$ 
et $k:X\times_S Y\to X\times Y$ les immersions ferm\'ees.
Les propri\'et\'es d'une famille de noyaux se sp\'ecialisent
\cite[Proposition 6.2]{Ch}:
\begin{prop}[Chen]
\label{famille}
Soit $K\in D^b(X\times_S Y)$ tel que $\Phi_{k_*K}:D^b(Y)\to D^b(X)$ est
pleinement fid\`ele (resp. est une \'equivalence).
Alors, $\Phi_{\BL(i\times j)^*K}:D^b(Y_0)\to D^b(X_0)$ est
pleinement fid\`ele (resp. est une \'equivalence).
\end{prop}

Ce r\'esultat permet de v\'erifier dans certains cas
qu'un noyau donne une \'equivalence par d\'eformation (cf par exemple
\S \ref{secflopMukai}).

\begin{rema}
Soit $\Phi_K:D^b(Y)\to D^b(X)$ un foncteur pleinement fid\`ele. Alors,
on doit penser \`a $Y$ comme l'espace de modules fin
de $\{F(\CO_y)\}_{y\in Y}$ et \`a $K$ comme l'objet universel associ\'e.
\end{rema}

\subsubsection{Repr\'esentabilit\'e}
\label{axiomes}
L'imperfection des axiomes des cat\'egories triangul\'ees rend la preuve du
r\'esultat suivant d\'elicate (cf \cite{OrK3}, \cite[Theorem 3.2.1]{Or2} et
\cite[Theorem 1.1]{BoVdB} qui assure l'existence d'adjoints;
cf \cite[Theorem 1.1]{Ka1} pour une extension aux champs de Deligne-Mumford
et une autre preuve).
\begin{theo}[Orlov]
\label{repOrlov}
Soit $F:D^b(Y)\to D^b(X)$ un foncteur pleinement fid\`ele. Alors,
il existe un unique $K\in D^b(X\times Y)$ tel que $F\simeq \Phi_K$.
\end{theo}

Une approche pr\'ef\'erable \`a ce probl\`eme (et \`a ceux de \S \ref{secinvariants})
consiste \`a consid\'erer une structure
plus riche que celle de cat\'egorie triangul\'ee, celle de dg-cat\'egorie
(les $\Hom$ sont munis d'une structure de complexe d'espaces vectoriels)
dont le \og $H^0$\fg\ est la cat\'egorie d\'eriv\'ee de d\'epart \cite{To}.
On dispose d'une dg-cat\'egorie
$L_{\coh}(X)$ dont le \og $H^0$\fg\ est $D^b(X)$.
Toen montre que la dg-cat\'egorie des
foncteurs de $L_{\coh}(Y)$ vers $L_{\coh}(X)$ est (quasi-)\'equivalente \`a 
$L_{\coh}(X\times Y)$ \cite[Theorem 8.15]{To} et ceci fournit un
analogue du th\'eor\`eme \ref{repOrlov}, pour des foncteurs non
n\'ecessairement pleinement fid\`eles. Les foncteurs $D^b(Y)\to D^b(X)$
obtenus sont alors tous du type $\Phi_K$ et r\'eciproquement tout foncteur
de ce type provient d'un foncteur d\'efini au niveau des dg-cat\'egories.

\subsection{Invariants d'une \'equivalence}
\label{secinvariants}
\subsubsection{}
Soit $F:D^b(Y)\iso D^b(X)$ une \'equivalence, avec $X$ et $Y$ projectives
lisses connexes. D'apr\`es le th\'eor\`eme \ref{repOrlov}, il existe
$K\in D^b(X\times Y)$ tel que $F\simeq\Phi_K$.

\medskip
La commutation de $F$ avec les foncteurs de Serre montre
que $\dim X=\dim Y$ et que $\omega_X$ et $\omega_Y$ ont le m\^eme ordre
\cite[Lemma 2.1]{BrMa}.

Un argument de rigidit\'e montre que $F$ induit un isomorphisme de
groupes alg\'ebriques
$\Pic^0(Y)\rtimes\Aut^0(Y)\iso \Pic^0(X)\rtimes\Aut^0(X)$,
o\`u $\Aut^0(X)$ est la composante neutre de $\Aut(X)$
(cf th\'eor\`eme \ref{equivabelienne} pour un cas o\`u les deux facteurs
sont \'echang\'es).

\subsubsection{Cohomologie}
\label{invariants}
Passons maintenant \`a des invariants du type cohomologie ou alg\`ebre
canonique (cf \cite[Theorem 4.9]{Mu2},
\cite{Ca1}, \cite{Ca2}, \cite[Theorem 2.3]{KaDK} et \cite[Theorem 2.1.8]{Or2}).

Soit
$HA_{i,k}(X)=\Ext^i_{X\times X}(\CO_{\Delta X},i_*\omega_X^k)$
o\`u $i:\Delta X\to X\times X$ est l'inclusion de la diagonale.
Soit $HA(X)=\bigoplus_{i,k}HA_{i,k}(X)$.
On munit $HA(X)$ d'une structure d'alg\`ebre bigradu\'ee via
les isomorphismes canoniques
$\Ext^i_{X\times X}(i_*\omega^r_X,i_*\omega^s_X)\iso
\Ext^i_{X\times X}(\CO_{\Delta X},i_*\omega^{s-r}_X)$.

On a $HA_{i,k}(X)=$\og$\Hom(S_X,S_X^k[i-k\dim X])$\fg, o\`u le terme
de droite doit \^etre compris comme le $H^0$ d'un complexe de $\Hom$'s
pris au niveau des dg-cat\'egories.

On a $HA_{i,k}(X)\simeq\bigoplus_{p+q=i}H^p(X,\Lambda^q\CT_X\otimes
\omega_X^k)$ (cf \cite{Ko2} et \cite[Corollary 2.6]{Sw}),
o\`u $\CT_X$ est le fibr\'e tangent.
En particulier, $\bigoplus_{k\ge 0}HA_{0,k}(X)$ est isomorphe \`a l'alg\`ebre
canonique $R(X)=\bigoplus_{k\ge 0}H^0(X,\omega_X^k)$.

\begin{theo}
\label{invHA}
$F$ induit un isomorphisme d'alg\`ebres gradu\'ees
$HA(Y)\iso HA(X)$.
En particulier, $F$ induit un isomorphisme gradu\'e entre les alg\`ebres canoniques
$R(Y)\iso R(X)$ et un isomorphisme entre les espaces vectoriels
de cohomologie $H^*(Y,\BC)\iso H^*(X,\BC)$.
\end{theo}

\noindent{\sc Preuve} (esquisse) ---
On utilise l'\'equivalence
$\Phi_{K\boxtimes L}:D^b(Y\times Y)\iso D^b(X\times X)$
o\`u 
$L=K^\vee\otimes p^*\omega_X[\dim X]\simeq
K^\vee\otimes q^*\omega_Y[\dim Y]$ vu dans $D^b(X\times Y)$.$\Box$

L'isomorphisme $H^*(X,\BC)\iso H^*(Y,\BC)$ n'est pas compatible \`a la
multiplication ni \`a la graduation
classique, en g\'en\'eral. Il est par contre compatible \`a la graduation
donn\'ee par ${^nH}(X,\BC)=\bigoplus_{p-q=n}H^p(X,\Omega^q_X)$.

\begin{rema}
\label{Hodge}
On s'attend tout de m\^eme \`a l'\'egalit\'e des nombres de Hodge de $X$
et $Y$ (cf \cite[\S 1.3]{BrMa2}). Dans la situation de la conjecture
\ref{conjflop}, la formule de changement de variable pour l'int\'egration
motivique montre que les vari\'et\'es $X$ et $Y$ ont les m\^emes nombres de
Hodge (Kontsevich, Batyrev et Denef--Loeser), cf \cite[\S 4]{Re} et \cite{Lo}.
\end{rema}

On d\'eduit du th\'eor\`eme \ref{invHA}
l'invariance de la dimension de Kodaira. On d\'eduit
aussi que si $\kappa(X,\omega_X)=\dim X$ (\ie, $X$ de type g\'en\'eral) ou
$\kappa(X,\omega_X^{-1})=\dim X$, alors $X$ et $Y$ sont birationnelles.

\medskip
Via le th\'eor\`eme de Grothendieck-Riemann-Roch, on obtient l'invariance de
la cohomologie \`a coefficients rationnels.

Soit $L\in D^b(X\times Y)$.
Soit $\ch:K_0(X)\to H^*(X,\BQ)$ la classe de Chern et soit
$\td_X$ la classe de Todd de $X$.
On d\'efinit un morphisme
$$\phi_L:H^*(Y,\BQ)\to H^*(X,\BQ),\ 
\xi\mapsto
 p_*\left(p^*(\sqrt{\td_X})\cdot \ch([L])\cdot q^*(\sqrt{\td_Y})\cdot
 q^*(\xi))\right).$$
Ce morphisme est gradu\'e pour la graduation ${^nH}$ et on a
$\phi_{L\circ L'}=\phi_L\circ\phi_{L'}$.

Lorsque $L=K$, alors ce morphisme est un isomorphisme et
sa complexification est l'isomorphisme du th\'eor\`eme \ref{invHA}.

La transformation $\Phi_L$ induit un morphisme
$[\Phi_L]:K_0(Y)\to K_0(X)$ et
on a un diagramme commutatif
$$\xymatrix{
K_0(Y)\ar[r]^{[\Phi_L]}\ar[d]_{\ch(-)\cdot\sqrt{\td_Y}} &
 K_0(X)\ar[d]^{\ch(-)\cdot\sqrt{\td_X}}\\
H^*(Y,\BQ)\ar[r]_{\phi_L} & H^*(X,\BQ)
}$$

\begin{rema}
Hille et Van den Bergh \cite[Remark 3.4]{HiVdB} mentionnent l'invariance
de la $K$-th\'eorie topologique par \'equivalence d\'eriv\'ee et en d\'eduisent
l'invariance de $H^*(X,\BZ)$ dans $H^*(X,\BQ)$ pour les courbes,
les surfaces $K3$ et les vari\'et\'es ab\'eliennes.
\end{rema}

\subsection{Vari\'et\'es ab\'eliennes}
Le r\'esultat suivant de Mukai \cite[Theorem 2.2]{Mu} est le point de d\'epart
des travaux sur les
\'equivalences entre cat\'egories d\'eriv\'ees de faisceaux coh\'erents.

Soient $A$ une vari\'et\'e ab\'elienne et $\hat{A}=\Pic^0(A)$ sa vari\'et\'e
ab\'elienne duale. On note $\CP$ le fibr\'e de Poincar\'e sur $A\times\hat{A}$.

\begin{theo}[Mukai]
\label{equivabelienne}
Le foncteur $\Phi_\CP:D^b(\hat{A}\mcoh)\to D^b(A\mcoh)$ est une \'equivalence.
\end{theo}

\noindent{\sc Preuve}  ---
On v\'erifie les conditions de la proposition \ref{critere}.
Pour $x\in \hat{A}$, alors
$\Phi_K(\CO_x)$ est un fibr\'e en droites $L_x$ de degr\'e $0$ sur $A$.
Puisque $H^*(L)=0$ pour tout $L\in \Pic^0(A)$ non trivial, on d\'eduit que
$\Phi_\CP$ est pleinement fid\`ele, et donc une \'equivalence,
car $\omega_A\simeq\CO_A$.$\Box$

On d\'ecrit maintenant toutes les \'equivalences d\'eriv\'ees entre vari\'et\'es
ab\'eliennes (cf \cite{Orab}, \cite[\S 5]{Or2} et \cite[\S 11 et \S15]{Po}).

Soit $B$ une vari\'et\'e ab\'elienne. 
Soit
$$f=\begin{pmatrix}x&y\\z&t\end{pmatrix}:B\times\hat{B}\to A\times\hat{A}.$$
On pose
$$\tilde{f}=\begin{pmatrix}\hat{t}&-\hat{y}\\-\hat{z}&\hat{x}
\end{pmatrix}:A\times\hat{A}\to B\times\hat{B}.$$
On note $U(B\times \hat{B},A\times\hat{A})$ l'ensemble des
isomorphismes $f:B\times\hat{B}\iso A\times\hat{A}$ tels que
$f^{-1}=\tilde{f}$.

Soient $a\in A$ et $\alpha\in\hat{A}$. On note
$m_a:A\to A,\ b\mapsto a+b$. On pose
$\Phi_{a,\alpha}=L_\alpha\otimes m_{a*}(-):D^b(A)\iso D^b(A)$.

\begin{theo}[Polishchuk, Orlov]
\label{equivab}
Soit $F:D^b(B)\iso D^b(A)$ une \'equivalence.
Alors, il existe $\sigma\in U(B\times \hat{B},A\times\hat{A})$ tel que
\begin{equation}
\label{swap}
\Phi_{\sigma(b,\beta)}=F\circ \Phi_{b,\beta}\circ F^{-1}
\textrm{ pour tous }b\in B \textrm{ et }\beta\in \hat{B}.
\end{equation}

R\'eciproquement, soit $\sigma\in U(B\times \hat{B},A\times\hat{A})$. Alors,
il existe une \'equivalence $F:D^b(B)\iso D^b(A)$ v\'erifiant (\ref{swap}).
\end{theo}

\noindent{\sc Preuve} (\'el\'ements) ---
L'invariance de $\Aut^0\times\Pic^0$ (cf \S \ref{secinvariants}) fournit
un isomorphisme $B\times\hat{B}\iso A\times\hat{A}$, dont on v\'erifie
qu'il a la propri\'et\'e voulue.
La r\'eciproque requiert la
construction d'un fibr\'e simple semi-homog\`ene de pente donn\'ee.$\Box$

On obtient alors une description explicite du groupe des auto-\'equivalences.
On a une suite exacte de groupes
$$0\to (A\times \hat{A})(\BC)\times\BZ \to \Aut(D^b(A))\to
U(A\times \hat{A},A\times \hat{A})\to 1.$$

\begin{rema}
On peut conjecturer qu'une vari\'et\'e projective lisse d\'eriv\'e-\'equivalente
\`a une vari\'et\'e ab\'elienne est une vari\'et\'e ab\'elienne.
\end{rema}

\subsection{Surfaces}

Notons tout d'abord que le cas des courbes n'est pas int\'eressant!
Deux courbes projectives lisses sont d\'eriv\'e-\'equivalentes si et seulement
si elles sont isomorphes~: en genre $\not=1$, le th\'eor\`eme
\ref{BO} donne le r\'esultat. Pour les courbes elliptiques, on le d\'eduit
par exemple de l'invariance de la structure de Hodge enti\`ere.

\smallskip
Nous d\'ecrivons bri\`evement la situation pour les surfaces.
Les d\'emonstrations demandent une analyse minutieuse suivant la
classification des surfaces minimales.

Soient $X$ et $Y$ deux surfaces projectives lisses connexes non isomorphes.
On suppose que $X$ est minimale, \ie, ne contient pas de $\BP^1$ avec
auto-intersection $-1$.
On a une description pr\'ecise des cas d'\'equivalences d\'eriv\'ees \cite{BrMa2}
(cf \cite{Mu2,OrK3} pour les $K3$).

\begin{theo}[Bridgeland-Maciocia]
On a $D^b(X)\simeq D^b(Y)$ si et seulement si une des assertions suivantes
est vraie
\begin{itemize}
\item $X$ et $Y$ sont toutes deux ab\'eliennes (ou toutes deux des $K3$)
et il existe une isom\'etrie entre
leurs r\'eseaux transcendants compatible avec les structures de Hodge
\item
$X$ et $Y$ sont des surfaces elliptiques et $Y$ est un sch\'ema de Picard
relatif de la fibration elliptique de $X$\cite{Brell}.
\end{itemize}
\end{theo}

\begin{rema}
\label{Uehara}
La $K$-\'equivalence entre deux vari\'et\'es projectives lisses $X$ et $Y$
(=existence d'un diagramme (\ref{Kequiv}) avec $f^*K_X\sim g^*K_Y$)
implique que les vari\'et\'es sont isomorphes en codimension $1$. Si deux
surfaces sont $K$-\'equivalentes, elles sont donc isomorphes.
Uehara \cite{Ue} construit des exemples de surfaces elliptiques
birationnelles et d\'eriv\'e-\'equivalentes mais qui ne sont pas isomorphes, donc 
pas $K$-\'equivalentes non plus.
\end{rema}

\section{Flips et flops}
\label{flipflop}
\subsection{Introduction}
Un cas particulier de la conjecture \ref{conjflop} est celui d'un flop.
Ce cas est important car une des conjectures du MMP est que deux
mod\`eles minimaux birationnels sont connect\'es par une suite de flops.

Un flop est un diagramme
$$\xymatrix{
X\ar[dr]_f & & X^+\ar[dl]^{f^+} \\
& \bar{X}
}$$
o\`u
\begin{itemize}
\item
$\bar{X}$ est une vari\'et\'e projective de Gorenstein,
\item
$f$ et $f^+$ sont des r\'esolutions cr\'epantes (\ie, $f^*\omega_{\bar{X}}\simeq
\omega_X$ et $(f^+)^*\omega_{\bar{X}}\simeq \omega_{X^+}$)
dont le lieu exceptionnel est de codimension $\ge 2$ et
\item
il existe un diviseur $D$ sur $X$ tel que $-D$ est
relativement $f$-ample et le transform\'e strict de $D$ est
relativement $f^+$-ample.
\end{itemize}
Le morphisme $f$ d\'etermine uniquement $f^+$ (car
$X^+=\Proj\bigoplus_{m\ge 0}\CO_{\bar{X}}(f_*(mD))$, ind\'ependant de $D$).
La conjecture \ref{conjflop} pr\'edit que $D^b(X)\simeq D^b(X^+)$ et
donc que
$X^+$ peut se construire comme un espace de module
d'objets de $D^b(X)$ (cf \S \ref{secdim3} pour la dimension $3$).
Cette approche pourrait aussi s'appliquer pour les flips.

\smallskip
Nous allons voir des exemples (\S \ref{secflopstandard},
\S \ref{goodMukai} et \S \ref{secdim3}) o\`u la transformation de noyau
$\CO_{X\times_{\bar{X}}X^+}$ est une \'equivalence $D^b(X^+)\iso D^b(X)$.
N\'eanmoins, nous verrons dans \S \ref{stratifie} une situation o\`u cette
transformation n'est pas une \'equivalence.

\subsection{\'Eclatements}
\subsubsection{Fibr\'es projectifs}
Soit $E$ un fibr\'e vectoriel
de rang $r\ge 1$ sur une vari\'et\'e projective lisse $Y$ et soit
$q:\BP(E)\to Y$ le fibr\'e projectif associ\'e.

La proposition suivante \cite[Theorem 2.6]{Or1} fournit une version relative
de l'exemple \ref{Beilinson}.
\begin{prop}[Orlov]
\label{fibreprojectif}
On a une d\'ecomposition semi-orthogonale
$$D^b(\BP(E))=\langle D^b(Y)_{-r+1},D^b(Y)_{-r+2},\ldots,D^b(Y)_{0}\rangle,$$
o\`u $D^b(Y)_d$ est l'image de $D^b(Y)$ par le foncteur pleinement
fid\`ele $\CO_q(d)\otimes q^*(-)$.
\end{prop}

\noindent{\sc Preuve} (esquisse) ---
La pleine fid\'elit\'e est fournie par la proposition
\ref{pleinefideliteprojection}.
Soient $C,D\in D^b(Y)$ et $L\in D^b(\BP(E))$. On a
$$\Hom(q^*C,L\otimes^\BL q^*D)\simeq
\Hom(C,\BR q_*(L\otimes^\BL q^*D))\simeq
\Hom(C,(\BR q_*L)\otimes^\BL D).$$
La semi-orthogonalit\'e r\'esulte alors de l'annulation de
$\BR q_*\CO_q(i)$ pour $-1\ge i\ge -r+1$.

Soit $\CT$ la sous-cat\'egorie triangul\'ee de $D^b(\BP(E))$
engendr\'ee par $ D^b(Y)_{-r+1},\ldots,D^b(Y)_{0}$.
Soient $y\in Y$ et $F=q^{-1}(y)$. Les faisceaux $\CO_F(-r+1),\ldots,
\CO_F(0)$ engendrent $D^b(F)$ (exemple \ref{Beilinson}), donc,
vus comme faisceaux sur $\BP(E)$, ils engendrent $D^b_F(\BP(E))$.
Par cons\'equent, $D^b_F(\BP(E))\subset \CT$ et en particulier
$\CO_x\in \CT$ pour $x\in F$. On en d\'eduit que ${^\perp\CT}=0$,
d'o\`u $\CT=D^b(\BP(E))$ par le th\'eor\`eme \ref{admissible}.$\Box$

\subsubsection{D\'ecomposition pour un \'eclatement}
Soient maintenant
$X$ une vari\'et\'e projective lisse et $Y$ une sous-vari\'et\'e ferm\'ee lisse
de $X$ purement de codimension $r\ge 2$.
Soit $\CN_{Y/X}$ le fibr\'e normal de $Y$ dans $X$.
Soient $\tX$ l'\'eclat\'ee de $X$ le long de $Y$ et 
$\tY\simeq \BP(\CN_{Y/X})$ le diviseur exceptionnel,
image inverse de $Y$ dans $\tX$.
On a un diagramme commutatif
$$\xymatrix{
\tY\ar@{^{(}->}[r]^j\ar[d]_q & \tX\ar[d]^p \\
Y\ar@{^{(}->}[r]_i & X 
}$$

Le th\'eor\`eme suivant d\'ecrit comment la cat\'egorie d\'eriv\'ee grossit par
\'eclatement \cite[Theorem 4.3]{Or1}.
\begin{theo}[Orlov]
\label{eclatement}
On a une d\'ecomposition semi-orthogonale
$$D^b(\tX)=\langle D^b(Y)_{-r+1},D^b(Y)_{-r+2},\ldots,
D^b(Y)_{-1},D^b(X)_0\rangle$$
o\`u $D^b(Y)_d=j_*(\CO_q(d)\otimes q^*D^b(Y))$ et $D^b(X)_0=\BL p^* D^b(X)$.
Les foncteurs canoniques induisent des \'equivalences
$D^b(X)\iso D^b(X)_0$ et
$D^b(Y)\iso D^b(Y)_d$ pour tout $d$.
\end{theo}

\noindent{\sc Preuve} (esquisse) ---
Un point utile est l'existence, pour tout $C\in D^b(\tY)$, d'un
triangle distingu\'e
\begin{equation}
\label{tri}
C\otimes \CO_q(1)\to \BL j^*j_* C\to C\rightsquigarrow.
\end{equation}

La proposition \ref{pleinefideliteprojection} montre
que $\BL p^*:D^b(X)\to D^b(\tX)$ est pleinement fid\`ele.
Pour l'\'equivalence $D^b(Y)\iso D^b(Y)_d$, la proposition
\ref{fibreprojectif} ram\`ene le probl\`eme \`a montrer que
$\Hom(\BL j^*j_* \CO_F,\CO_F[i])=0$ pour toute fibre $F$ de
$q$ et tout $i\not=0$. Le triangle distingu\'e (\ref{tri})
ram\`ene l'annulation recherch\'ee \`a celles de $H^*(\BP^{r-1},\CO(-1))$ et
$H^{>0}(\BP^{r-1},\CO)$.

\smallskip
La semi-orthogonalit\'e de 
$\CT=\langle D^b(Y)_{-r+1},D^b(Y)_{-r+2},\ldots, D^b(Y)_{-1},D^b(X)_0\rangle$
s'\'etablit par la m\^eme technique que dans la preuve de la
proposition \ref{fibreprojectif}.

Il reste \`a \'etablir que ${^\perp \CT}$ est nul.
Soit $y\in Y$. On a $\BL^dp^*\CO_{y}\simeq \Omega_{p^{-1}(y)}^d(d)$ pour
$0\le d\le r-1$ et $\BL^dp^*\CO_{y}=0$ pour $d\ge r$. On a
$H^*(\BP^{r-1},\Omega^d(d))=0$ pour $1\le d\le r-1$ et
la r\'esolution de la diagonale dans $\BP^{r-1}$ (cf l'exemple \ref{Beilinson})
montre que pour $1\le d\le r-1$, alors $\Omega^d(d)$ est dans la
sous-cat\'egorie de $D^b(\BP^{r-1})$ engendr\'ee par
$\CO(-r+1),\ldots,\CO(-1)$. Par cons\'equent, le c\^one du
morphisme canonique $\BL p^*\CO_y\to \CO_{p^{-1}(y)}$ est dans $\CT$, donc
finalement $\CO_{p^{-1}(y)}$ est dans $\CT$ et on termine comme pour
la proposition \ref{fibreprojectif}.$\Box$

\begin{rema}
Les arguments essentiels pour les deux r\'esultats pr\'ec\'edents apparaissent
aussi dans les travaux de Thomason \cite{Thfib,Thecl} qui s'int\'eresse
\`a la $K$-th\'eorie sup\'erieure.
\end{rema}

\subsection{Flips et flops standard}
\label{secflopstandard}

Soient $X$ une vari\'et\'e projective lisse et $Y$ une sous-vari\'et\'e ferm\'ee isomorphe
\`a $\BP^k$ telle que $\CN_{Y/X}\simeq \CO_Y(-1)^{l+1}$ avec $1\le l\le k$.
Soit $p:\tX\to X$ l'\'eclat\'ee de $X$ le long de $Y$. Le diviseur exceptionnel
$\tY$ est isomorphe \`a $\BP^k\times\BP^l$ et on a $\CN_{\tY/\tX}\simeq
\CO(-1,-1)$. Soit
$p^+:\tX\to X^+$ la contraction de $\tY$ sur $\BP^l$. On suppose que
$X^+$ est une vari\'et\'e projective.
Soient $f:X\to \bar{X}$ et
$f^+:X^+\to \bar{X}$ les contractions de $\BP^k$ et $\BP^l$ sur
un point. Ils fournissent l'exemple le plus simple de flip (et de flop
lorsque $k=l$).
On a $\tX\simeq X\times_{\bar{X}}X^+$.
$$\xymatrix{
 & \tY\simeq \BP^k\times\BP^l\ar@{^{(}->}[d]^j\ar[dr]\ar[dl]\\
\BP^k\simeq Y\ar@{^{(}->}[d] & \tX\ar[dl]^-p \ar[dr]_-{p^+} &
 Y^+\simeq\BP^l\ar@{^{(}->}[d]\\
X\ar[dr]_f && X^+\ar[dl]^{f^+} \\
& \bar{X}
}$$

Les conjectures \ref{conjflop} et
\ref{conjflip} sont connues dans ce cas
(cf \cite[Theorem 3.6]{BoOr1} et \cite[Theorem 2.2.9]{Or2}):

\begin{theo}[Bondal-Orlov]
Le foncteur $\Phi_{\CO_\tX}=\BR p_*\BL(p^+)^*:D^b(X^+)\to D^b(X)$ est
pleinement fid\`ele et c'est une \'equivalence si $l=k$.
\end{theo}

\noindent{\sc Preuve} (esquisse) ---
Soit $\CA$ (resp. $\CB$) la sous-cat\'egorie triangul\'ee pleine de $D^b(\tX)$
engendr\'ee par les $j_*\CO(r,s)$ avec $r=0$ et
$-l\le s\le -1$ (resp.  $-k\le r<0$ et $-l\le s\le -1$).
Les d\'ecompositions du th\'eor\`eme \ref{eclatement} et de l'exemple
\ref{Beilinson} donnent une d\'ecomposition semi-orthogonale
$(\BL p^*D(X))^\perp=\langle\CB,\CA\rangle$.
On a
$\CB\subset(\BL(p^+)^*D^b(X^+))^\perp$. Puisque
$(\omega_\tX)_{|\tY}\simeq \CO(-k,-l)$, on a aussi
$\CA\otimes\omega_\tX\subset(\BL(p^+)^*D^b(X^+))^\perp$, donc
$\CA\subset{^\perp(\BL(p^+)^*D^b(X^+))}$.

Soit $C\in D^b(X^+)$ et soit $\bar{C}$ le c\^one du morphisme canonique
$\BL p^*\BR p_*(\BL(p^+)^* C)\to \BL(p^+)^* C$. On a
$\bar{C}\in (\BL p^*D(X))^\perp\cap{^\perp\CB}=\CA$.
Par cons\'equent, le morphisme canonique
$$\Hom(C,D)\to \Hom(\BR p_* \BL (p^+)^* C,\BR p_* \BL(p^+)^* D)$$
est un isomorphisme pour tout $D\in D^b(X^+)$.

Supposons maintenant $k=l$. Soit $p^!=\BL p^*(-)\otimes\omega_{\tX/X}[k]$.
Le foncteur
$\BR p_*^+ p^!$ est adjoint \`a droite de $\BR p_* \BL(p^+)^*$.
On montre, comme ci-dessus, que
$\BR p_*^+ (\omega_{\tX/X}[k]\otimes\BL p^*(-))$ est pleinement fid\`ele et
finalement $\BR p_* \BL(p^+)^*$ est une \'equivalence.$\Box$

\begin{rema}
Il devrait en fait \^etre vrai que $D^b(X)/\BR p_*\BL (p^+)^*D^b(X^+)$
a une suite exceptionnelle compl\`ete de longueur $k-l$.
\end{rema}

\subsection{Flop de Mukai}
\label{secflopMukai}
\subsubsection{}
\label{goodMukai}
Soit $X$ une vari\'et\'e projective lisse de dimension $2n\ge 4$. Soit
$i:Y=\BP^n\hookrightarrow X$ une immersion ferm\'ee.
On suppose $\CN_{Y/X}\simeq \Omega^1_Y$. On consid\`ere comme
pr\'ec\'edemment l'\'eclatement $p:\tX\to X$. Soit
$Y^\vee$ l'espace projectif dual. Alors,
le diviseur exceptionnel $\tY$
s'identifie \`a la vari\'et\'e d'incidence dans $Y\times Y^\vee$ et on a
$\CN_{\tY/\tX}\simeq \CO(-1,-1)_{|\tY}$.
On a alors une contraction $p^+:\tX\to X^+$ de $\tY$ sur $Y^\vee$ et
on suppose $X^+$ projective.

Le foncteur $\BR p_*\BL(p^+)^*:D^b(X^+)\to D^b(X)$ n'est pas
une \'equivalence (cf \cite[Proposition 5.12]{KaDK} et
\cite[Corollary 2.2]{Na1}).
Par contre, un autre noyau fournit une
\'equivalence (cf \cite[Corollary 5.7]{KaDK} et \cite[Theorem 4.4]{Na1}).

Soient $f:X\to \bar{X}$ et
$f^+:X^+\to \bar{X}$ les contractions de $Y$ et $Y^+$ sur
un point (flop de Mukai).
Soit $\hat{X}=X\times_{\bar{X}}X^+$. Alors, 
$\hat{X}=\tX\cup(Y\times Y^+)\subset X\times X^+$ et l'intersection est \`a
croisements normaux.

$$\xymatrix{
& \hat{X}\ar[dl]_\pi \ar[dr]^{\pi^+} \\
X \ar[dr]_-f && X^+\ar[dl]^-{f^+} \\
& \bar{X}
}$$

\begin{theo}[Kawamata, Namikawa]
\label{flopMukai}
Le foncteur
$\Phi_{\CO_{\hat{X}}}=\BR\pi_* \BL(\pi^+)^*:D^b(X^+)\to D^b(X)$ est une
\'equivalence.
\end{theo}

\noindent{\sc Preuve} (esquisse) ---
Il suffit de traiter le cas o\`u $X=\Spec\Omega^1_Y$ et $i:Y\to X$
est la section nulle, car cela ne change pas le compl\'et\'e formel
de $X$ le long de $Y$ (cf \S \ref{secplf}). La non-projectivit\'e de $X$
ne pose pas de probl\`eme nouveau. Soit
$\CX=\Spec\CO_Y(-1)^{n+1}$.
La suite exacte $0\to \Omega_Y^1\to \CO_Y(-1)^{n+1}\to \CO_Y\to 0$
fournit un morphisme lisse $\CX\to\BA^1$ et la fibre de $0$ s'identifie
\`a $X$.

On a $\CN_{Y/\CX}\simeq\CO_Y(-1)^{n+1}$ et on a alors un flop
standard $\tilde{\CX}\to\CX^+$ (cf \S \ref{secflopstandard}).
Le foncteur $\Phi_{\CO_{\tilde{\CX}}}:D^b(\CX^+)\to D^b(\CX)$
est une \'equivalence et il r\'esulte de la proposition \ref{famille} que
$\Phi_{\CO_{\hat{X}}}$ est une \'equivalence.$\Box$

Deux vari\'et\'es projectives symplectiques birationnelles de dimension $4$
sont connect\'ees par une suite de flops de Mukai 
(\cite[Theorem 1.2]{WieWis}, \cite[\S 8]{CMS}; 
il faut en fait
permettre la contraction simultan\'ee de plusieurs $\BP^2$ disjoints
pour que les vari\'et\'es interm\'ediaires restent projectives
\cite{Wie}).
On d\'eduit alors (\cite[Corollary 4.5]{Na1} et \cite[Remark 5.13]{KaDK})~:

\begin{coro}[Kawamata, Namikawa]
\label{sym4}
Deux vari\'et\'es projectives symplectiques birationnelles de dimension $4$
sont d\'eriv\'e-\'equivalentes.
\end{coro}

\subsubsection{Flop de Mukai stratifi\'e}
\label{stratifie}

Markman \cite{Ma} a \'etudi\'e une g\'en\'eralisation du flop de Mukai.
Nous suivons ici la pr\'esentation de \cite{Na2} et ne donnons que
la version \og lin\'earis\'ee\fg.

Soit $V$ un espace vectoriel de dimension $n$ et soit
$G(V,r)$ la Grassmannienne des sous-espaces de dimension $r$ de $V$,
o\`u $r$ est un entier $\le n/2$.
Soient $X=T^*G(V,r)$ et $\bar{X}$ la sous-vari\'et\'e ferm\'ee de
$\End_{\BC}(V)$ des endomorphismes $a$ tels que $a^2=0$ et $\rang(a)\le r$.
On identifie $T^*G(V,r)$ \`a la vari\'et\'e des paires
$(W,\phi)$ o\`u $W$ est un sous-espace de dimension $r$ de $V$ et
$\phi\in\Hom(V/W,W)$.
On dispose d'une application moment
$f:X\to \bar{X}$ qui envoie $(W,\phi)$
sur la composition
$a:V\twoheadrightarrow V/W\xrightarrow{\phi}W\hookrightarrow V$.
L'application moment est un isomorphisme au-dessus de l'ouvert des
endomorphismes de rang $r$.
On dispose de m\^eme d'une application moment
$f^+:X^+=T^*G(V,n-r)\to \bar{X}$.

Le noyau \og \'evident\fg\ ne fournit pas une \'equivalence d\'eriv\'ee,
bien que sa classe dans $K_0$ soit ad\'equate \cite[Theorem 2.7 et
Observation 4.9]{Na2}:

\begin{theo}[Namikawa]
L'application $[\Phi_{\CO_{X\times_{\bar{X}}X^+}}]:K_0(X^+)\to K_0(X)$ est un
isomorphisme pour tous $n,r$ mais
$\Phi_{\CO_{X\times_{\bar{X}}X^+}}:D^b(X^+)\to D^b(X)$ n'est pas une
\'equivalence pour $n=4$ et $r=2$.
\end{theo}

Dans le cas $n=4$ et $r=2$, Kawamata construit un noyau de
la forme $\CO_{X\times_{\bar{X}}X^+}\otimes\CL$ o\`u $\CL$ est un fibr\'e en
droites, qui induit une \'equivalence \cite{KaG}.

\smallskip
Voyons maintenant une construction pour $n,r$ g\'en\'eraux d'un
noyau support\'e par $X\times_{\bar{X}}X^+$
et induisant une \'equivalence d\'eriv\'ee.

Soit $U$ la sous-vari\'et\'e ouverte de $G(V,r)\times G(V,n-r)$ des
$(W,W')$ tels que $W\cap W'=0$ et soit
$\iota:U\to G(V,r)\times G(V,n-r)$ l'immersion ouverte.
$$
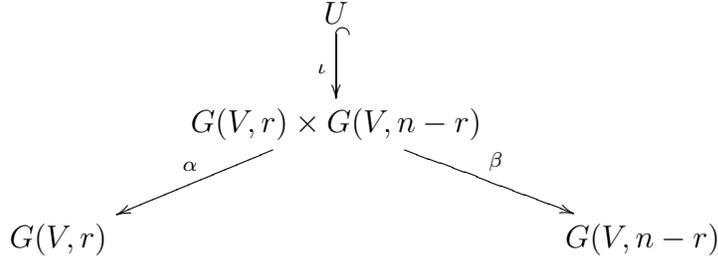
\xymatrix{
& U\ar@{^{(}->}[d]_\iota \\
& G(V,r)\times G(V,n-r)\ar[dl]_{\alpha}\ar[dr]^{\beta} \\
G(V,r) && G(V,n-r)
}$$
Il est classique \cite[Exercice III.15]{KaScha} que la transformation
de noyau $\iota_!\BC_U$ induit une \'equivalence entre cat\'egories
d\'eriv\'ees de faisceaux constructibles de $\BC$-espaces vectoriels:
$$\BR\alpha_!(\BC_U\otimes \beta^{-1}(-)):
D^b_{\textrm{cons}}(G(V,n-r))\iso D^b_{\textrm{cons}}(G(V,r))$$
Soit $\CK$ le module de Hodge mixte correspondant \`a $\iota_!\BC_U$.
C'est un noyau inversible pour
les transformations entre cat\'egories d\'eriv\'ees de modules de Hodge mixtes.
Soit $K=\mathrm{Gr}(\CK)$, un faisceau coh\'erent
$\BG_m$-\'equivariant sur $T^*(G(V,r)\times G(V,n-r))=X\times X^+$. Alors,
l'inversibilit\'e du noyau $\CK$ montre que
$\Phi_K:D^b(X^+)\to D^b(X)$ est une \'equivalence (c'est bien s\^ur
aussi une \'equivalence pour les cat\'egories $\BG_m$-\'equivariantes).

Il serait int\'eressant de d\'ecrire explicitement le faisceau $K$. 
Kashiwara a propos\'e une description pour $n=4$ et $r=2$, qu'il
faudrait comparer avec le noyau de \cite{KaG}.

Cette construction se g\'en\'eralise \`a des vari\'et\'es de drapeaux paraboliques
pour des groupes semi-simples complexes arbitraires. Dans le cas
des vari\'et\'es de drapeaux complets, une telle construction au niveau
de la $K$-th\'eorie a \'et\'e effectu\'ee par Tanisaki \cite{Ta}.

On peut aussi construire une \'equivalence $D^b(X^+)\iso D^b(X)$ en
regroupant les $D^b(T^*G(V,r))$, $0\le r\le \dim V$, et en
appliquant les m\'ethodes de \cite{ChRou}.

\subsection{Dimension $3$}
\label{secdim3}
\subsubsection{}
On expose ici la construction de Bridgeland
\cite{Brdim3} (voir aussi \cite{Ch,KaDK,VdB}).

Consid\'erons un flop entre
vari\'et\'es projectives lisses de dimension $3$~:
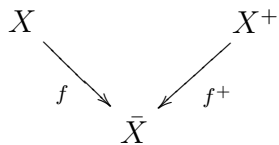
$$\xymatrix{
X\ar[dr]_f & & X^+\ar[dl]^{f^+} \\
& \bar{X}
}$$

\begin{rema}
Lorsque $f$ ne contracte qu'une courbe irr\'eductible $C$, alors
$C\simeq\BP^1$ et le fibr\'e normal $\CN_{C/X}$
est $\CO(-1)\oplus\CO(-1)$, $\CO\oplus\CO(-2)$ ou $\CO(1)\oplus\CO(-3)$
(cf \cite[Corollary 16.3]{ClKoMo}). Le premier cas est un flop
standard (\S \ref{secflopstandard}) et le second cas peut se traiter par
des m\'ethodes similaires \cite[Theorem 3.9]{BoOr1}.
\end{rema}

Consid\'erons le diagramme
$$\xymatrix{
& X\times_{\bar{X}} X^+\ar[dl]_\pi \ar[dr]^{\pi^+} \\
X && X^+
}$$

\begin{theo}[Bridgeland]
\label{equivdim3}
Le foncteur
$\BR \pi_* \BL (\pi^+)^*:D^b(X^+)\to D^b(X)$ est une \'equivalence.
\end{theo}

Bridgeland construit le flop $f^+$ \`a partir de $f$~:
la vari\'et\'e $X^+$ appara{\^\i}t comme un espace de modules fin de certains objets de
$D^b(X)$ (\og faisceaux pervers ponctuels\fg) et le noyau de
l'\'equivalence est le fibr\'e universel. La d\'etermination de ce fibr\'e, et
donc la forme pr\'ecise du th\'eor\`eme \ref{equivdim3}, sont dues \`a Chen.

Un flop g\'en\'eralis\'e entre vari\'et\'es projectives lisses de dimension $3$ se
d\'ecompose en une suite de flops \cite[Theorem 4.6]{KaDK} et
on en d\'eduit:

\begin{coro}
\label{CY3}
La conjecture \ref{conjflop} est vraie en dimension $3$.
\end{coro}

\subsubsection{Construction du flop et \'equivalence}
Soit $\bar{X}$ une vari\'et\'e projective connexe de Gorenstein de dimension $3$.
Soit $f:X\to\bar{X}$ une r\'esolution cr\'epante dont le lieu exceptionnel est
union d'un nombre fini de courbes et soit $D$ un diviseur sur $X$ tel que
$-D$ est relativement $f$-ample.

Le th\'eor\`eme d'annulation de Grauert-Riemenschneider montre, via la formule
de projection, que $\BR^{>0} f_*\CO_X=0$.
Le foncteur $\BL f^*:D(\bar{X})\to D(X)$ est donc pleinement fid\`ele
(variante de la proposition \ref{pleinefideliteprojection} pour les
cat\'egories d\'eriv\'ees non born\'ees que requiert la non-lissit\'e de $\bar{X}$).
Soit $\CB$ son image.
On a une d\'ecomposition semi-orthogonale
$D(X)=\langle \CB^\perp,\CB\rangle$, o\`u $\CB^\perp=\{C\in D(X)|\BR f_* C=0\}$.

On construit une nouvelle t-structure sur $D(X)$ en recollant la
t-structure standard de $\CB$ avec celle de $\CB^\perp$ d\'ecal\'ee de
$1$ vers la gauche.
Son c{\oe}ur $\Per_{X/\bar{X}}$ (une cat\'egorie ab\'elienne)
consiste en les $C\in D(X)$ tels que
\begin{itemize}
\item $\CH^iC=0$ pour $i\not=0,1$
\item $f_*\CH^{-1}C=0$
\item $\BR^1f_*\CH^0C=0$ et $\Hom(\CH^0C,M)=0$ pour tout
$M\in \CB^\perp\cap X\mcoh$.
\end{itemize}

On a $\CO_X\in \Per_{X/\bar{X}}$ et on dit que $E\in\Per_{X/\bar{X}}$ est un
\og faisceau pervers ponctuel\fg\ si
\begin{itemize}
\item $E$ est un quotient de $\CO_X$ (dans $\Per_{X/\bar{X}}$)
\item la classe de Chern de $E$ est \'egale \`a celle d'un faisceau
gratte-ciel.
\end{itemize}

Soit $S$ une vari\'et\'e. Une
famille de faisceaux pervers ponctuels param\'etr\'ee par $S$ est
un objet $\CE\in D^b(S\times X)$ tel que $\BL j^*_s\CF\in\Per_{X/\bar{X}}$
pour tout point (ferm\'e) $s\in S$. Ici, $j_s:s\times X\to S\times X$
est l'inclusion. Deux familles qui diff\`erent par le produit par
un fibr\'e inversible de $S$ sont dites \'equivalentes.

On consid\`ere le foncteur qui \`a une vari\'et\'e $S$ associe l'ensemble
des classes d'\'equivalence de familles de faisceaux pervers ponctuels
param\'etr\'ees par $S$. Bridgeland \cite[Theorem 3.8]{Brdim3} d\'emontre
l'existence d'un espace de modules fin:
le foncteur est repr\'esentable par une vari\'et\'e projective $M(X/\bar{X})$.

Soit $U$ l'ouvert de $\bar{X}$ des points au-dessus desquels $f$ est un
isomorphisme.
Alors, $\CO_{f^{-1}(u)}$ est un faisceau pervers ponctuel pour tout $u\in U$.
Ceci d\'efinit une immersion ouverte $U\to M(X/\bar{X})$ et on note
$X^+$ l'adh\'erence de l'image (on montre en fait que $X^+=M(X/\bar{X})$).
On note $f^+:X^+\to \bar{X}$ le morphisme canonique.
Chen \cite[Proposition 4.2]{Ch} montre que
$\CO_{X\times_{\bar{X}} X^+}$ est un fibr\'e universel.

\smallskip
Le th\'eor\`eme \ref{equivdim3} admet la version plus pr\'ecise suivante.
\begin{theo}[Bridgeland]
\label{BKR}
Le foncteur
$\BR \pi_* \BL (\pi^+)^*:D^b(X^+)\to D^b(X)$ est une \'equivalence.
En outre, $f^+:X^+\to \bar{X}$ est un flop~: $f^+$ est
une r\'esolution cr\'epante et le transform\'e strict
de $D$ est relativement $f^+$-ample.
\end{theo}

La preuve est essentiellement la m\^eme que pour la
correspondance de McKay \cite{BKR}.
L'outil-clef est un r\'esultat d'alg\`ebre commutative que nous rappelons
maintenant (dans la version de \cite[\S 5]{BrMa2}).

Soient $Z$ une vari\'et\'e irr\'eductible et $C\in D^b(Z)$ non nul. Soit
$\Supp(C)$ l'union des supports des $\CH^i(C)$.
Soit $\ampl(C)$ l'ensemble des $i\in\BZ$ tels qu'il existe $z\in Z$
avec $\Hom(C,\CO_z[-i])\not=0$.
La dimension homologique de $C$ est $\hd(C)=\sup(\ampl(C))-\inf(\ampl(C))$.

\begin{theo}[\og Nouveau th\'eor\`eme d'intersection\fg]
\label{nouveau}
On a $\codim \Supp(C)\le \hd(C)$.

Soit $z\in Z$ tel que $\Supp(C)=\{z\}$, $\CH^0(C)\simeq\CO_z$ et
$\ampl(C)\subseteq [-\dim Z,0]$. Alors, $z$ est un point lisse de $Z$ et
$C\simeq \CO_z$.
\end{theo}

\noindent{\sc Preuve du th\'eor\`eme \ref{BKR}} (esquisse) ---
Soit $p_X:X\times X^+\to X$ la premi\`ere projection.
Soient $\CP=\CO_{X\times_{\bar{X}} X^+}$ et $\CP_w=\Phi_\CP(\CO_w)$ pour
$w\in X^+$. Soient $\CP'=\CP^\vee\otimes p_X^*\omega_X[3]$ et
$\CQ=\CP'\circ\CP\in D^b(X^+\times X^+)$.
On v\'erifie que $\Hom(\CP_w,\CP_{w'}[i])=0$ pour tout $i$,
lorsque $f^+(w)\not=f^+(w')$
et on montre que $\Hom(\CP_w,\CP_{w'})=\delta_{w,w'}\BC$, d'o\`u
$\Hom(\CP_{w'},\CP_w[3])\simeq\delta_{w,w'}\BC$ par dualit\'e de Serre.

Soit $\CQ'=\CQ_{|X^+\times X^+-\Delta X^+}$. Le support de $\CQ'$ est
contenu dans $X^+\times_{\bar{X}} X^+-\Delta X^+$, donc est de codimension
$\ge 2$, s'il est non vide. On a $\Hom(\CQ',\CO_{w,w'}[i])\simeq
\Hom(\CP_w,\CP_{w'}[i])$, donc
$\hd(\CQ')\le 1$. Le
th\'eor\`eme \ref{nouveau} montre alors que $\CQ'=0$.

Le morphisme canonique
$\Hom(\CO_w,\CO_w[1])\to
\Hom(\CP_w,\CP_w[1])$
est injectif (il
s'identifie \`a l'application de Kodaira-Spencer). On en d\'eduit que
$\CH^0(\Phi_{\CQ}(\CO_w))\simeq\CO_w$. On a
$\ampl(\Phi_{\CQ}(\CO_w))=[-3,0]$. Le 
th\'eor\`eme \ref{nouveau} montre que
$\Phi_{\CQ}(\CO_w)\simeq\CO_w$ et que $X^+$ est lisse en $w$.
La cr\'epance de $f$ permet alors de conclure que $\Phi_\CP$ est 
une \'equivalence, via la proposition \ref{critere}.

La cr\'epance de $f^+$ s'obtient en utilisant la trivialit\'e du foncteur de
Serre de $D^b_F(X^+)$ pour toute fibre $F$ de $f^+$ \cite[Lemma 3.1]{BKR}.
On montre que $\Phi_\CP[-1]$ se restreint en une \'equivalence entre
faisceaux coh\'erents sur $X^+$ tu\'es par $\BR f^+_*$ et
faisceaux coh\'erents sur $X$ tu\'es par $\BR f_*$ et 
l'amplitude relative du transform\'e strict de $D$ s'en d\'eduit
\cite[\S 4.6]{Brdim3}.$\Box$

\begin{rema}
Il serait int\'eressant de voir si la m\^eme m\'ethode fournit une
construction de flips (cf \cite{AbCh1} pour des r\'esultats dans cette
direction, en dimension quelconque).
\end{rema}

\subsection{Vari\'et\'es singuli\`eres}

Les constructions pour les vari\'et\'es projectives lisses ont des
g\'en\'eralisations \`a une classe plus large, dans le cadre du MMP
\cite{AbCh2,Ch,KaDK,Kaflip,Ka1,VdB}.
Ceci est utile pour construire des \'equivalences entre vari\'et\'es projectives
lisses, car il existe des d\'ecompositions de flops g\'en\'eralis\'es en flops qui
font intervenir des vari\'et\'es singuli\`eres.
Kawamata d\'emontre en particulier
un th\'eor\`eme d'\'equivalence d\'eriv\'ee
pour des flops particuliers entre champs de Deligne-Mumford
toriques. Il en d\'eduit qu'un flop torique g\'en\'eralis\'e donne lieu
\`a une \'equivalence d\'eriv\'ee \cite[Corollary 4.5]{Ka1}:
\begin{theo}[Kawamata]
\label{torique}
Soient $f:Z\to X$ et $g:Z\to Y$ des morphismes toriques birationnels entre
vari\'et\'es toriques projectives lisses tels que $f^*\omega_X\simeq
g^*\omega_Y$. Alors,
$D^b(Y)\simeq D^b(X)$.
\end{theo}

Mentionnons pour terminer la n\'ecessit\'e de consid\'erer des vari\'et\'es
analytiques et d'effectuer les constructions dans le cadre de la 
log-g\'eom\'etrie.

\end{document}